\documentclass[preprint,12pt,authoryear]{elsarticle}

\usepackage{amssymb}

\usepackage{lineno}
\usepackage{amsmath, tabu}
\usepackage[usenames, dvipsnames]{color}
\usepackage{amstext} 
\usepackage{array}   
\usepackage{graphicx}
\usepackage{caption}
\usepackage{subcaption}
\usepackage{algorithm}
\usepackage{algorithmic}

\journal{Journal}

\begin{document}

\begin{frontmatter}



\title{Nonlinear Expectation Inference for Efficient Uncertainty Quantification and History Matching of Transient Darcy Flows in Porous Media with Random Parameters Under Distribution Uncertainty}

\address[1]{Research Centre for Mathematics and Interdisciplinary Sciences, Shandong University, Qingdao, Shandong Province, 266237, China}
\address[2]{Frontiers Science Center for Nonlinear Expectations, Minister of Education, Shandong University, Qingdao, Shandong Province, 266237, China}
\address[3]{Institute of Mathematical Sciences, ShanghaiTech University, Pudong, Shanghai, 201210, China}
\address[4]{School of Petroleum Engineering, China University of Petroleum (East China), Qingdao, Shandong Province, 266580, China}
\address[5]{School of Civil Engineering, Qingdao University of Technology, Qingdao, Shandong Province, 266520, China}

\author[1,2]{Zhao Zhang\corref{cor1}}
\cortext[cor1]{Corresponding author}
\ead{zhaozhang@sdu.edu.cn}
\author[1,2]{Xinpeng Li}
\author[1,2]{Menghan Li}
\author[3]{Jiayu Zhai}
\author[5]{Piyang Liu}
\author[4]{Xia Yan}
\author[4,5]{Kai Zhang}

\begin{abstract}
The uncertainty quantification of Darcy flows using history matching is important for the evaluation and prediction of subsurface reservoir performance. Conventional methods aim to obtain the maximum a posterior or maximum likelihood estimate (MLE) using gradient-based, heuristic or ensemble-based methods. These methods can be computationally expensive for high-dimensional problems since forward simulation needs to be run iteratively as physical parameters are updated. In the current study, we propose a nonlinear expectation inference (NEI) method for efficient history matching and uncertainty quantification accounting for distribution or Knightian uncertainty. Forward simulation runs are conducted on prior realisations once, and then a range of expectations are computed in the data space based on subsets of prior realisations with no repetitive forward runs required. In NEI, no prior probability distribution for data is assumed. Instead, the probability distribution is assumed to be uncertain with prior and posterior uncertainty quantified by nonlinear expectations. The inferred result of NEI is the posterior subsets on which the expected flow rates are consistent with observation. The accuracy and efficiency of the new method are validated using single- and two-phase Darcy flows in 2D and 3D heterogeneous reservoirs.  

\end{abstract}


\end{frontmatter}


\section{Introduction}
\label{intro}
The uncertainty quantification (UQ) of Darcy flows in porous media is important in many physical and engineering problems involving groundwater flows, geo-energy reservoirs and CO2 storage \citep{oliver2008inverse, caers2011modeling, smith2013uncertainty, feng2024deep}. Parameters of subsurface Darcy flows, e.g. permeability and porosity, can be highly heterogeneous, but measurements of these parameters are only available at limited locations. In consequence, these parameters are highly uncertain in geological models.

Geostatistical methods have been developed to estimate the distributions of uncertain parameters and generate a range of prior realisations to account for uncertainty \citep{soares2001direct, boyd2019quantifying, strebelle2020multiple}. Geostatistical methods can be grouped into two-point and multi-point methods. Two-point methods quantify spatial correlations using variogram models and generate realisations stochastically. The main disadvantage of two-point methods is the difficulty to reproduce complex geological structures, e.g. channels. For such scenarios, multi-point methods can be adopted to generate geologically meaningful realisations using training images that honour the geological features.

Simulation results using prior realisations are often different from observed data, and it is hence necessary to calibrate geological models such that simulation results can match observed dynamic responses. This calibration process is called history matching in reservoir engineering and is a prerequisite for making predictions based on simulation results \citep{oliver2008inverse, caers2011modeling, arnold2019uncertainty}. Generally, there are manual and automatic history matching methods. Manual methods adjust parameters manually and are very time consuming. Automatic history matching (AHM) calibrate parameters using optimisation algorithms which are generally more efficient than manual methods \citep{zhang2016two}. 

AHM methods typically aim to obtain the maximum a posterior (MAP) or maximum likelihood estimate (MLE) for physical parameters. Gradient-based, heuristic and ensemble-based methods can be used for optimising the parameters. The computation of Jacobian and Hessian matrix in gradient-based methods, e.g. Newton, quasi-Newton and steepest descent, is very time consuming \citep{yan2021stein}. Heuristic methods such as evolution approaches and simulated annealing requires no gradients but the convergence rate is often very slow \citep{ma2020multiscale}. Ensemble-based methods, e.g. ensemble Kalman filter (EnKF) and ensemble smoother multiple-data assimilation (ESMDA), optimise physical parameters using cross-covariance matrices, thus being much more efficient and can quantify the uncertainty in dynamic responses using posterior realisations \citep{emerick2013ensemble, schillings2023ensemble}. But for high-dimensional problems, ensemble-based methods need more realisations and correspondingly more forward simulation runs. Further, Markov Chain Monte Carlo and variational Bayes inference could be adopted to approximate the posterior probability density, but they can be computationally expensive for high-dimensional problems \citep{cui2016dimension, olalotiti2018multiobjective, povala2022variational}.

A major difficulty of AHM is that expensive forward simulation needs to be conducted iteratively as physical parameters are optimised \citep{zhang2021history}. To reduce the computational costs of forward simulation, data-driven or physics-informed surrogates can be built to predict dynamic responses efficiently given physical parameters as inputs \citep{newsum2017efficient, zhang2022physics, ma2022novel, zhang2023physics,guan2023fast,yan2024physics}. In addition, inversion algorithms based on the Bayes' rule such as EnKF and ESMDA requires prior knowledge about the probability distribution of geological parameters.

In the current study, we propose a new nonlinear expectation inference (NEI) method for the UQ and history matching of transient Darcy flows with random parameters, guided by the nonlinear expectation theory \citep{peng2019nonlinear,peng2023improving}. An advantage of NEI for practical applications is that NEI assumes the probability distribution of parameters is uncertain, referred to as 'distribution uncertainty' or 'knightian uncertainty'. 
The notion of distribution uncertainty is more general than the usual case where we assumes a certain probability distribution function with a few uncertain hyperparameters, e.g. mean and variance. Distribution uncertainty acknowledge the fact that the form of probability distribution function itself is uncertain.
This is practical for the UQ of geological parameters as their probability distributions are highly uncertain given limited observation and experiments. 
Distribution uncertainty has been acknowledged and studied in many previous works, mostly in economics and finance \citep{nishimura2007irreversible, rindova2020shape, townsend2024chance}. Nonlinear expectation can be adopted as a robust approach to quantity distribution uncertainty, as evidenced by recent studies \citep{peng2019law,cohen2019data,li2022maximally,ji2023imbalanced,chen2023strategic}. 

Apart from the robustness regarding probability distribution, NEI is also very efficient requiring forward simulations only on a set of prior realisations while the inference is done in the data space using expectations built on subsets of realisations. The inferred result is subsets of realisations whose expected responses approximately cover the varying interval of observation. This implies that posterior sample distributions, i.e. subsets, are the results of inference, and the expected responses on subsets are to be used as the 'correct' prediction. This is different from most existing methods where a single expectation is obtained using responses on posterior realisations.

The paper is organised as follows. In section 2, the uncertainty quantification problem for subsurface flows is introduced. In section 3, the new nonlinear expectation inference (NEI) method is developed. In section 4, NEI is validated on 2D and 3D test cases of transient Darcy flows. Section 5 is the conclusion of the paper.

\section{The Uncertainty Quantification Problem}
\subsection{Forward Models}
In the current study, we focus on the UQ of random parameters associated with PDEs governing physical processes. 
Let $m\in \mathcal{M}$ be the random parameter and $d\in \mathcal{D}$ be the observations of solution, where $\mathcal{M}$ and $\mathcal{D}$ are Banach spaces, we have in compact form \citep{povala2022variational}

\begin{equation}
	d=\mathcal{G}(m)+\epsilon~,\label{first}
\end{equation}
where $\mathcal{G}$ is the mapping from $\mathcal{M}$ to $\mathcal{D}$, and $\epsilon$ is additive observation noise. Let $\mathcal{A}: \mathcal{M}\rightarrow\mathcal{U}$ be a solution operator of the governing PDE, where $\mathcal{U}$ is a Hilbert space, and $\mathcal{P}: \mathcal{U}\rightarrow\mathcal{D}$ be a projection operator, Eq.~\eqref{first} can be written as

\begin{equation}
	d=\mathcal{P}(\mathcal{A}(m))+\epsilon~.\label{second}
\end{equation}

\subsubsection{Single-Phase Transient Darcy Flow}
The governing equation for transient single-phase slightly compressible Darcy flows is \citep{chen2006book} 

\begin{equation}
\phi c_t\frac{\partial p}{\partial t}=\nabla \cdot(\frac{K}{\mu}\nabla p)+q~, \label{pt}
\end{equation}
where $p=p(\textbf{x}, t)$ is pressure, $K=K(\textbf{x})$ is permeability, $c_t$ is total compressibility, $\mu$ is viscosity, $\phi=\phi(\textbf{x})$ is porosity, $q=q(\textbf{x})$ is the source/sink term, $t$ is time and $\textbf{x}$ is the spatial coordinate. The finite volume method with implicit time integration can be adopted for solving Eq.~\eqref{pt} numerically as
\begin{equation}
V_i\phi C_t\frac{p_i^{n+1}-p_i^n}{\Delta t}=\sum_jT_{ij}(p_j^{n+1}-p_i^{n+1})+q_i^{n}V_i ~,\label{dis}
\end{equation}
where $V_i$ is the volume of the $i$'th grid cell, $p_i^n$ is the pressure at the $i$'th cell and time step $n$. $T_{ij}$ is the transmissibility between two neighbouring cells evaluated by harmonic average as

\begin{equation}
T_{ij}=(T_i^{-1}+T_j^{-1})^{-1}~,
\end{equation}
where $T_i$ is the transmissibility inside cell $i$ towards $j$ calculated as

\begin{equation}
T_i=\frac{K_iA_{ij}}{\mu d_i}~,
\end{equation}
where $A_{ij}$ is the area of the boundary face between cells $i$ and $j$, while $d_i$ is the distance from the cell centre to the face centre.
The volumetric flow rate $q_i^{n}V$ is evaluated at cells with non-zero sink term by

\begin{equation}
q_i^{n}V_i=\text{PI}*(p_i^{n}-p_{wf}) \label{wellmodel}
\end{equation}
where $p_{wf}$ is the constant bottom-hole pressure and PI is the production index \citep{peaceman1978}.

\subsubsection{Two-Phase Transient Darcy Flow}
The governing equations for transient two-phase slightly compressible Darcy flows are, neglecting gravity and capillary pressure \citep{chen2006book},
	
\begin{equation}
\phi\left[S_\alpha(c_r+c_\alpha)\frac{\partial p}{\partial t}+\frac{\partial S_\alpha}{\partial t}\right]=\nabla\left(\frac{k_{r\alpha}K}{\mu_\alpha}\nabla p_\alpha \right)=q_\alpha ~,
\end{equation}
where $\alpha$ is $o$ for the non-wetting and $w$ for wetting phase, $S_\alpha$ is phase saturation, $k_{r\alpha}$ is phase relative permeability, $c_r$ and $c_\alpha$ are rock and phase compressibility, respectively. The finite volume method and implicit-pressure explicit-saturation time integration scheme are adopted for discretization. Summing up the equations for wetting and non-wetting phases yields the governing equation for pressure
	
\begin{equation}
\phi(c_r+S_oc_o+S_wc_w)\frac{\partial p}{\partial t}=\nabla\left(\left(\frac{k_{ro}}{\mu_o}+\frac{k_{rw}}{\mu_w}\right)K\nabla p\right)=q_o+q_w~,
\end{equation}
which is discretized  as
	
\begin{equation}
\frac{V_i\phi c_t}{\Delta t}(p_i^{n+1}-p_i^n)=\sum_{j}\lambda^n_{ij} T_{ij}(p_j^{n+1}-p_i^{n+1})+(q_o+q_w)_i^nV_i~,
\end{equation}
where $\lambda^n_{ij}=(\lambda_o+\lambda_w)^n_{ij}=\left(\frac{k_{ro}}{\mu_o}+\frac{k_{rw}}{\mu_w}\right)^n_{ij}$ is the total mobility evaluated by the first-order upwind scheme at time step $n$, while $T_{ij}=\left(\frac{KA}{d}\right)_{ij}$ is computed by harmonic average. Saturation and flow rates are fixed for computing pressure. Then saturation is updated explicitly as
	
\begin{equation}
\frac{V_i\phi S_{w, i}^n C_w}{\Delta t}(p_i^{n+1}-p_i^n)+V\phi\frac{S^{n+1}_{w, i}-S^n_{w, i}}{\Delta t}=\sum_{j}\lambda^n_{w, ij} T_{ij}(p_j^{n+1}-p_i^{n+1})+q_wV_i ~.
\end{equation}
Then, phase mobilities and flow rates are updated for the next time step. The relative permeabilities are calculated using the Brooks-Corey model assuming zero residual oil saturation \citep{brooks1964}
	
\begin{eqnarray}
k_{rw}&=&\left(\frac{S_w-S_{iw}}{1-S_{iw}}\right)^{\frac{2+3\beta}{\beta}}~,\label{krw}\\
k_{ro}&=&\left(\frac{1-S_w}{1-S_{iw}}\right)^2\left(1-\left(\frac{S_w-S_{iw}}{1-S_{iw}} \right)^{\frac{2+\beta}{\beta}} \right)~,\label{kro}
\end{eqnarray}
where $\beta$ indicates the distribution of pore sizes that $\beta>2$ is for narrow distributions and $\beta<2$ for wide distributions. Here we set $\beta=2$. 

\subsection{Bayesian Uncertainty Quantification}
Although NEI is not built on the Bayes' rule, we introduce Bayesian UQ here for comparison and as a background for UQ methods. The infinite Bayesian approach is to build the inference model in the function space before discretizing the model, e.g. the pre-conditioned Crank-Nicholson MCMC scheme \citep{pinski2015algorithms}. However, the more adopted Bayesian approach is to build the inference model in the finite-dimensional space after discretizing the forward model \citep{deveney2023deep}, where the random parameter $m$ can be estimated using Bayes' rule as

\begin{equation}
	p(m|d)=\frac{p(d|m)p(m)}{\int p(d|m)p(m)dm}\propto p(d|m)p(m) ~,\label{bayes}
\end{equation}
where $p(m|d)$ is the posterior, $p(d|m)$ is the likelihood and $p(m)$ is the prior pdf. The likelihood is evaluated using the forward model given the pdf of observation noise. 

For nonlinear inverse problems, analytical solutions of the posterior pdf are generally not available and iterative inference methods need to be used.
Existing iterative inference algorithms are either to find a point-estimate or to approximate the posterior pdf itself. The former mainly includes gradient-based, heuristic and data assimilation methods for computing the maximum a posterior (MAP) or maximum likelihood estimate (MLE), while the latter includes MCMC and variational Bayes approaches \citep{posselt2013markov, emerick2013ensemble, povala2022variational, schillings2023ensemble}.
 
Bayesian methods for nonlinear inverse problems typically need repetitive evaluations of the forward model in the inference process, which can be computationally expensive for high-dimensional problems. To approximate the posterior pdf for UQ, prior pdf and likelihood function (e.g. Gaussian) needs to be known which is often not satisfied in practical problems. 

\section{Nonlinear Expectation Inference}
In this section, we develop the nonlinear expectation inference (NEI) method accounting for distribution uncertainty. The development of NEI is enlightened by the $\varphi$-max-mean algorithm in the nonlinear expectation theory \citep{peng2019nonlinear, ji2023imbalanced, peng2023improving}. 
According to the nonlinear expectation theory, the stochastic process for the physical parameter is considered to be of complex nature that no single distribution can serve as a perfect model, but should be described by a family of distributions.
Distribution uncertainty is a practical assumption for heterogeneous geological models where only sparse measurements of physical parameters are available, as evidenced by the studies using Gaussian mixture model to estimate the distribution of geophysical parameters \citep{li2017adaptive, astic2021petrophysically, wang2022gaussian}. 

Given no explicit form of prior distributions in practical problems, the uncertainty of distribution can be quantified using nonlinear expectations \citep{peng2019nonlinear}. Typical examples of nonlinear expectation can be defined as the upper and lower limits over a range of linear expectations \citep{gl}. Let $\mathcal{P}$ be a set of distributions characterizing distribution uncertainty, the upper expectation $\hat{E}$ is defined as the upper limit of linear expectations $E_P[X]$ over $\mathcal{P}$, i.e.,
$\hat{E}[X]=\sup_{P\in\mathcal{P}}E_P[X]$. Let $\psi$ be any given function, we have $\hat{E}[\psi(X)]=\sup_{P\in\mathcal{P}}E_P[\psi(X)]$.
Such upper expectation $\hat{E}$ is a sublinear expectation and linear additivity is no longer true for $\hat{E}$ as $\hat{E}[\psi(X)+\psi(Y)]\leq \hat{E}[\psi(X)]+\hat{E}[\psi(Y)]$ for random variables $X$ and $Y$ under distribution uncertainty. The lower expectation is superlinear accordingly but can be defined using sublinear expectation as $-\hat{E}[-\psi(X)]$. Hereafter, 'linear' expectation is referred to as expectation for simplicity.

Precisely, we consider the inverse problem of inferring $m$ given observed $d$ associated with the forward model in Eq.~\eqref{second} rewritten as

\begin{equation}
	d=\mathcal{P}(\mathcal{A}(m))+\epsilon=\mathcal{\psi}(m)+\epsilon~.\label{psi}
\end{equation}
where $\psi$ can be approximated by numerical simulation. The conventional assumption in AHM is that the probability distribution of $\epsilon$ is known while that of $m$ is unique but with unknown parameters. Here, we assume that the distribution of $m$ is not unique and aim to infer $m$ given $d$ and a number of prior realisations. 
The NEI method is developed as follows.

\noindent\textbf{(1) Estimate the Prior Uncertainty}

Let $\mathcal{M}=\{m_i, 1\leq i\leq n\}$ denote the set containing $n$ prior realisations. Enlightened by the $\psi$-max-mean algorithm \citep{peng2019nonlinear}, a range of expectations can be computed over subsets of $\mathcal{M}$ to characterise distribution uncertainty. The total number of subsets of $\mathcal{M}$ is $2^n$. Since $n>>1$, $2^n>>n$ and we have an enormous number of subsets. A subset is regarded as a possible sample distribution for $m$ such that the expectations on subsets characterise distribution uncertainty. The expectation on subset $M_j$ is denoted as $E_j$ for simplicity.

Practically, evaluating expectations on a total number of $2^n$ subsets may be too expensive, and we can use a limited number of subsets for estimation.
Suppose a number of $C$ subsets are built, we have $C$ expectations accordingly. The prior distribution uncertainty of $m$ can be quantified using the upper and lower expectations $\bar{\psi}=\hat{E}[\psi(m)]$ and $\underline{\psi}=-\hat{E}[-\psi(m)]$ such that

\begin{equation}
	\bar{\psi}=\max_{j=1,\dots, C}E_{M_j}[\psi(m)]=\max_{j=1,\dots, C}E_j, \label{upper}
\end{equation}

and 
\begin{equation}
	\underline{\psi}=-\hat{E}[-\psi(m)]=\min_{j=1,\dots, C}E_j, \label{lower}
\end{equation}

The true $\psi(m)$ is bounded by $\bar{\psi}$ and $\underline{\psi}$ given sufficient prior realisations.
Two random variables $X$ and $Y$ under distribution uncertainty are identical if their upper and lower expectations are identical, respectively.

Effectively, $\bar{\psi}$ and $\underline{\psi}$ can be obtained simply by taking the upper and lower limits of all $\psi_m$ over all realisations considered to be subsets containing only one element. But the computation of $E_j$ over subsets containing more than one elements is necessary to approximate the posterior uncertainty.

\noindent\textbf{(2) Estimate the Posterior Uncertainty}

Given observed data $d_{obs}$ and a loss function $L$, the loss of a realisation $m_i$ is evaluated using $L(d_{sim}, d_{obs})$ where $d_{sim}=\psi(m_i)$. The loss for each expectation on subset $M_j$  is 

\begin{equation}
	L_j=L(E_j, d_{obs})~.
\end{equation}

The posterior expectations $E^*_j$ and subsets $M^*_j$  are obtained by searching in all prior expectations $E_j$ for those satisfying a given criteria, e.g. $L_j<\sigma$ with $\sigma$ being a threshold adjustable according to the uncertainty of observed data. Here, $*$ marks posterior. The inferred result of NEI is the posterior subsets on which the expected responses are consistent with the varying interval of observation. 

The posterior upper and lower expectations quantifying distribution uncertainty are

\begin{eqnarray}
	\bar{\psi}^*=\max_{M^*_j} L^*_j& \le & \min_{M_j} L_j~, \nonumber\\
	\underline{\psi}^*=\min_{M^*_j} L^*_j& \ge & \min_{M_j} L_j~.\label{biglike}
\end{eqnarray}

The key for NEI is that by building a large number of prior expectations based on limited prior realisations, we can find a range of posterior expectations that satisfy the given criteria. 
The algorithm of NEI is summarised in Algorithm \ref{alg1}. The NEI method is highly efficient as we only need to conduct forward simulation runs on prior realisations for once, while the evaluation of linear expectations on subsets are in the data space with no further forward runs. For prediction with NEI, forward model runs only need to be conducted on non-repetitive realisations in the union of posterior subsets to compute expectations on posterior subsets accordingly.

The prerequisite for NEI is that the forward model is approximately correct, which is also the basis for other history matching methods, e.g. gradient-based and data assimilation methods. Another prerequisite for NEI is that the varying interval of observation be covered by responses simulated using prior realisations. In cases where this prerequisite is not fulfilled, more prior realisations need to be generated and simulated. In practise, it is also difficult for many other history matching methods (e.g. ESMDA and particle filtering) to converge if the observed data falls outside the range of the responses using prior realisations.

\begin{algorithm}
	\caption{Nonlinear Expectation Inference}
	\label{alg1}
	\begin{algorithmic}[1]
		\REQUIRE Prior realisation set $\mathcal{M}=\{m_i, 1\leq i\leq n\}$, observed data vector $d$, the loss function $L$ and the forward model $\psi$.
		\ENSURE Posterior linear, nonlinear expectations in the data space and associated subsets of realisations in the parameter space.
		\STATE Run forward model on all prior realisations to obtain the dynamic response $\psi(m)$ for each realisation.
		\STATE Build a number of $C$ subsets of $\mathcal{M}$. 
		\STATE Calculate the expectations of $\psi(m)$ on the subsets as $E_{M_j}[\psi(m)]=\frac{1}{k}\sum_{i=1}^k\psi(m_{ji})$, where $m_{ji}$ is the $i'$th realisation in the $j'$th subset.
		\STATE Calculate the prior upper $\bar{\psi}=\hat{E}[\psi(m)]$ and lower expectation $\underline{\psi}=-\hat{E}[-\psi(m)]$. The observation data $d$ should be bounded by $\bar{\psi}$ and $\underline{\psi}$. If not, introduce more prior realisations until the condition is true.  
		\STATE Evaluate the loss for the expectation associated with each subset.
		\STATE Obtain a number of $C^*$ posterior expectations and their related subsets, where $C^*$ is the number of posterior subsets with $L<\sigma$ where $\sigma$ is the threshold.
        \STATE Compute the upper and lower expectations of the $C^*$ posterior expectations, determine whether they fit within the upper and lower bounds of the observed data, and accordingly select an appropriate $\sigma$.
		\STATE The posterior upper and lower expectations are obtained as limits of posterior expectations to quantify distribution uncertainty. 
		\STATE For prediction, the non-repetitive realisations in the union of posterior subsets are simulated further to compute linear expectations on posterior subsets.
	\end{algorithmic}
\end{algorithm}

\section{Test Cases}
In the test cases, the forward model is the governing equation for transient Darcy flows. The uncertain parameter is the heterogeneous permeability field, and the observed data is a vector of flow rates for all time steps at the producing wells. In the first two test cases, the random parameter field is close to Gaussian for validation where the distribution uncertainty is mainly about hyperparameters. In the third test case, the probability distribution of the random parameter field is highly uncertain.

\subsection{Test Case 1: Single-Phase Transient Flow in 2D Heterogeneous Reservoir}
The first example is the inversion of a 2D 100 m$\times$100 m permeability field $k=k(x, y)$ given the observed flow rate (sink term). The permeability field is close to Gaussian for comparison with the ensemble smoother multiple data assimilation (ESMDA) \citep{emerick2013ensemble} method. ESMDA is a widely adopted efficient approach for history matching.

\begin{figure}[!htp] 
	\centering
	\begin{subfigure}[t]{0.49\textwidth}
		\centering
		\includegraphics[width=1\textwidth]{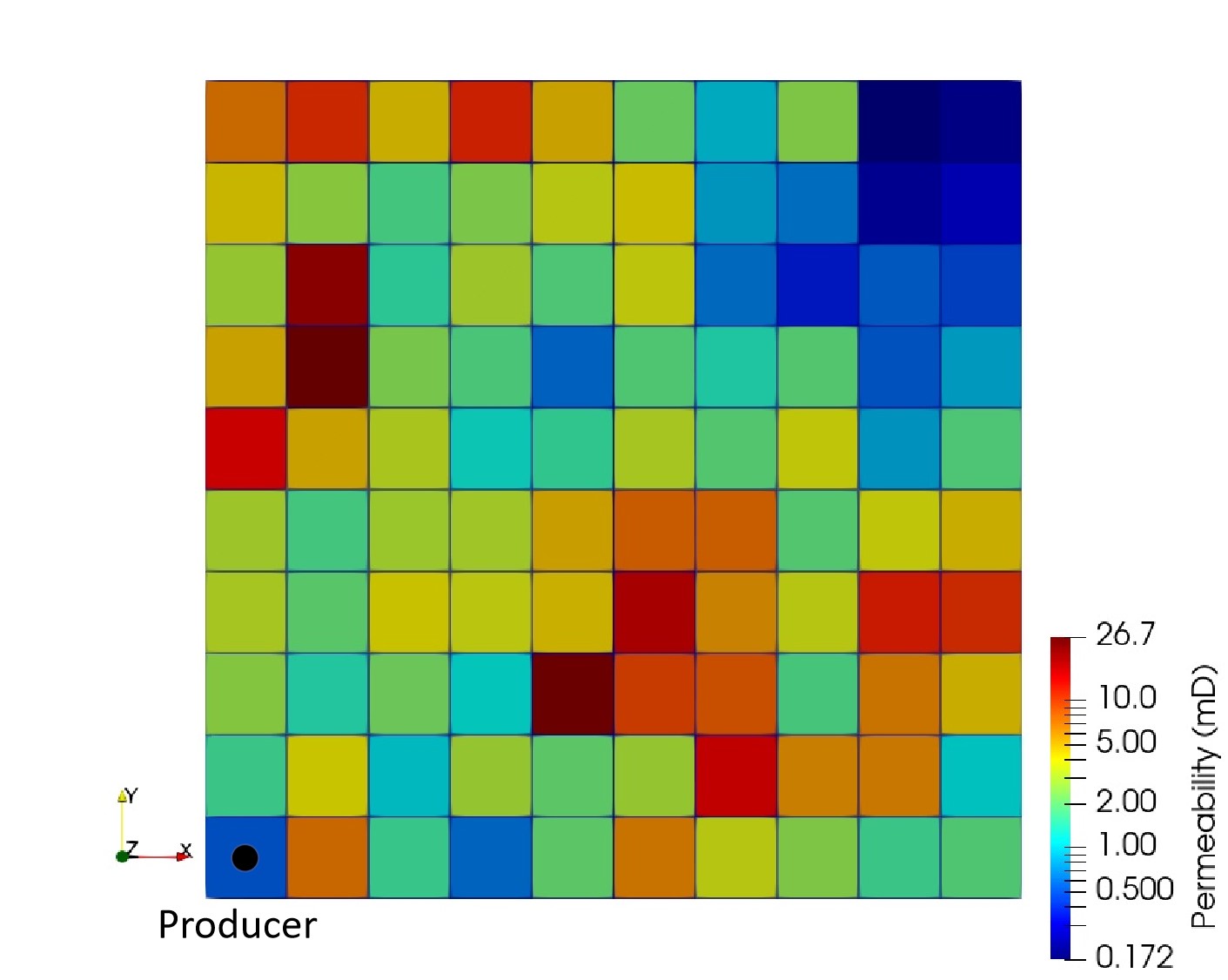}
		\caption{}
		\label{trueperm} 
	\end{subfigure}
	\begin{subfigure}[t]{0.49\textwidth}
		\centering
		\includegraphics[width=1\textwidth]{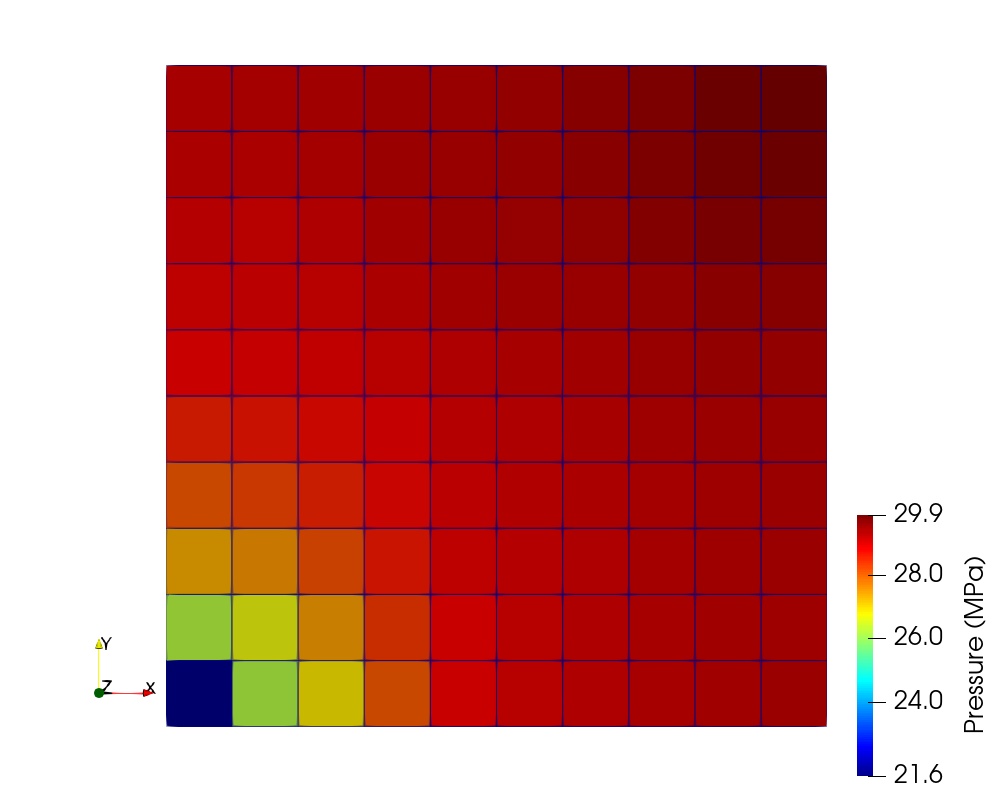}
		\caption{}
		\label{truepress}
	\end{subfigure}
	\caption{The ground truth for generating the observed flow rate data (a), and the pressure field at end of simulation using the true permeability field (b) for test case 1.}
\end{figure}

The grid is 10$\times$10 for validation. The permeability field is heterogeneous, isotropic and uncertain. All other parameters are constant shown in Table \ref{para1}. A total of 50 prior realisations for the permeability field are built via conducting sequential Gaussian simulation (SGS) using the Stanford Geostatstical Modeling Software (SGeMS) \citep{remy2009applied} for $h=log_2(10k)$ where the unit of permeability $k$ is mD. The ground truth of the permeability field for generating the synthetic observation data is presented in Fig.~\ref{trueperm}. All boundaries are with no-flow condition while the flow is driven by the transient sink term at a corner. The pressure field at end of simulation using the true permeability field is shown in Fig.~\ref{truepress}. For forward simulation, there are 180 time steps with $dt=100$k seconds. 

For NEI, we conduct forward simulation on the 50 prior realisations, while no further forward runs are needed in the inversion process. Let each subset contain no more than three realisations, a total of $\sum_{i=1}^{3}C_{50}^i=20875$ subsets are built to compute the expectations $E_i(\psi(k))$, $i=1,...,20875$ for flow rates $q=\psi(k)$. Expectations on subsets with less realisations are computed first, as they cover a wider range and the evaluation of less realisations is more efficient.
Next, set the threshold $\sigma$ for loss and calculate the posterior upper and lower expectations to determine whether they fit the upper and lower bounds of the observed data. Based on this, select an appropriate $\sigma$. 
The inferred result of NEI is the posterior subsets on which the expected flow rates cover the varying interval of observation. For prediction, the nonrepetitive realisations in posterior subsets are simulated further and the corresponding expectations are computed. 

For ESMDA, $h$ is inverted to recover $k$ according to $h=log_2(10k)$, as the distribution of $h$ is close to Gaussian. The ensemble of 50 realisations are updated for six iterations until convergence to yield the posterior realisations. For prediction using ESMDA, all 50 posterior realisations are simulated further.

\begin{figure}[!htp] 
	\centering
	\begin{subfigure}[t]{0.48\textwidth}
		\centering
		\includegraphics[width=1\textwidth]{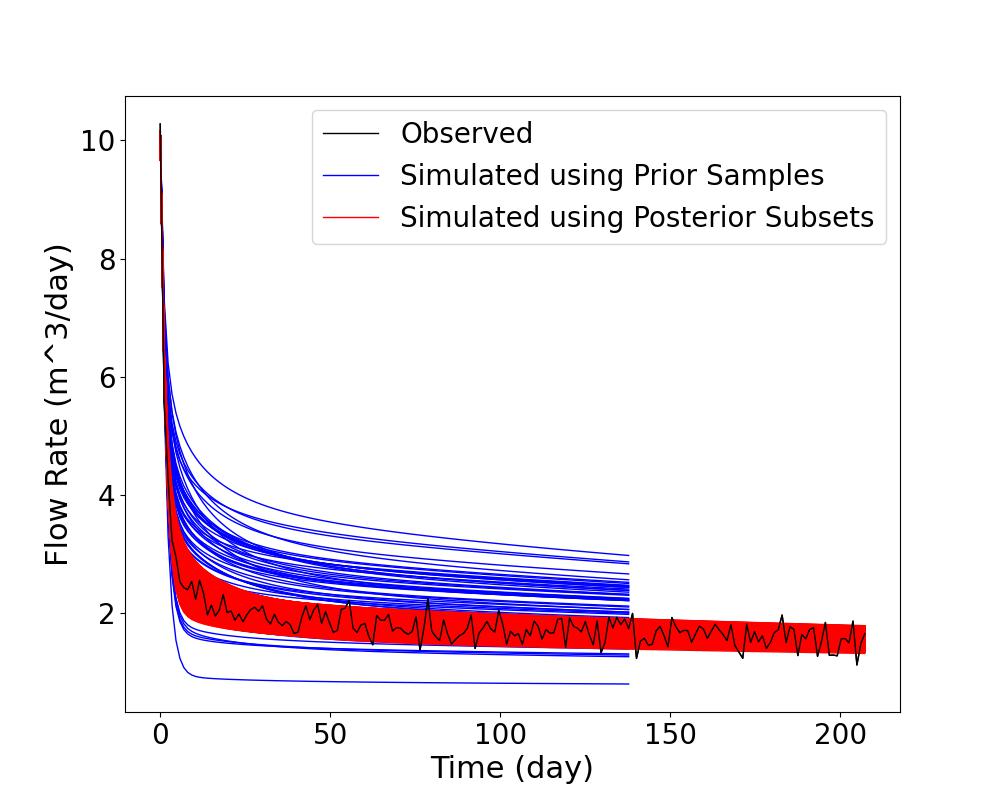}
		\caption{}
		\label{comparepriorpost} 
	\end{subfigure}
	\begin{subfigure}[t]{0.48\textwidth}
		\centering
		\includegraphics[width=1\textwidth]{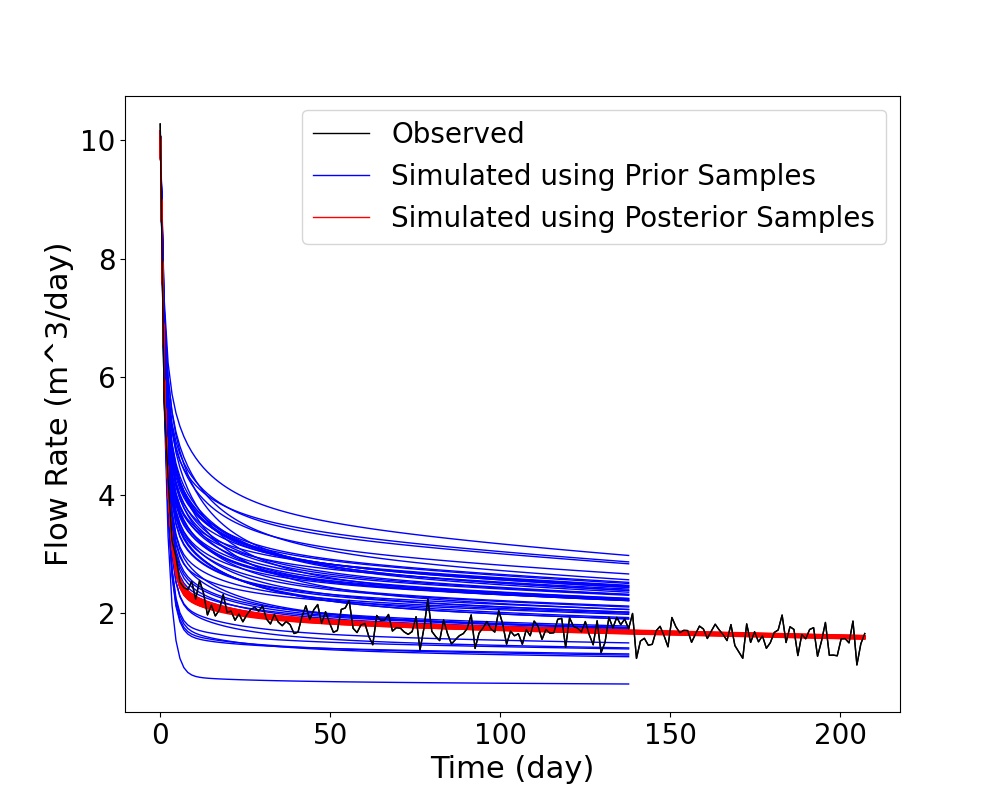}
		\caption{}
		\label{comparepriorpost2}
	\end{subfigure}
	\caption{The predicted flow rate at wells using posterior subsets by NEI (a) and the predicted flow rate at wells using posterior realisations by ESMDA (b) for test case 1.}
	\label{flow2d}
\end{figure}

The inference results of NEI and ESMDA can be seen from Fig.~\ref{flow2d} presenting the flow rates simulated using prior and posterior subsets or realisations. The posterior flow rates of NEI is the expectations on posterior subsets.
In addition, NEI is computationally more efficient using 50 forward simulation runs while ESMDA needs an extra 300 forward runs. 
The cpu time for ESMDA is about 26 seconds, and that for NEI is about 6 seconds, where 2.5 seconds is for forward simulation and 3.5 seconds is for computing the expectations on subsets (Table \ref{time}). The computation of all test cases are on a desktop workstation with i9 central processing unit.

\begin{table}[!htp]
	\centering
	\begin{tabular}{|l|l|}
		\hline
		Porosity & 10\%\\
		\hline
		Total Compressibility & $5\times 10^{-8}$ Pa$^{-1}$\\
		\hline
		Dynamic Viscosity & 0.002 Pa*s\\
		\hline
		Production Index & 1.175$\times$10$^{-5}$ m$^3$Sec$^{-1}$MPa$^{-1}$\\
		\hline
		Bottom-hole Pressure & 20 MPa \\
		\hline
	\end{tabular}
	\caption{A summary of constant physical and well parameters for test case 1.}
	\label{para1}
\end{table}

\begin{table}[!htp]
	\centering
	\begin{tabular}{|l|l|l|l|l|}
		\hline
		Test cases & Dimension & Phases & NEI cost & ESMDA cost \\
		\hline
		Test Case 1 & 100 & Single & 6 sec & 26 sec\\
		\hline
		Test Case 2 & 2000 & Two & 6.4 hours & 13.5 hours \\
		\hline
		Test Case 3 & 25200 & Single & 16 hours & doesn't converge\\
		\hline
	\end{tabular}
	\caption{CPU time of NEI and ESMDA on test cases.}
	\label{time}
\end{table}

\subsection{Test Case 2: Two-Phase Transient Flow in 3D Heterogeneous Reservoir}
The second example is the inversion of a 3D 300 m$\times$300 m$\times$30 m permeability field $k=k(x, y, z)$ given the observed flow rates at producing wells for two-phase Darcy flow. The grid is 20$\times$20$\times$5. A total of 100 prior realisations are generated by SGS for $h=log_2(10k)$ where the unit of permeability $k$ is mD. The permeability field is close to Gaussian to compare NEI and ESMDA.

\begin{figure}[!htp]
	\centering
	\includegraphics[width=0.6\linewidth]{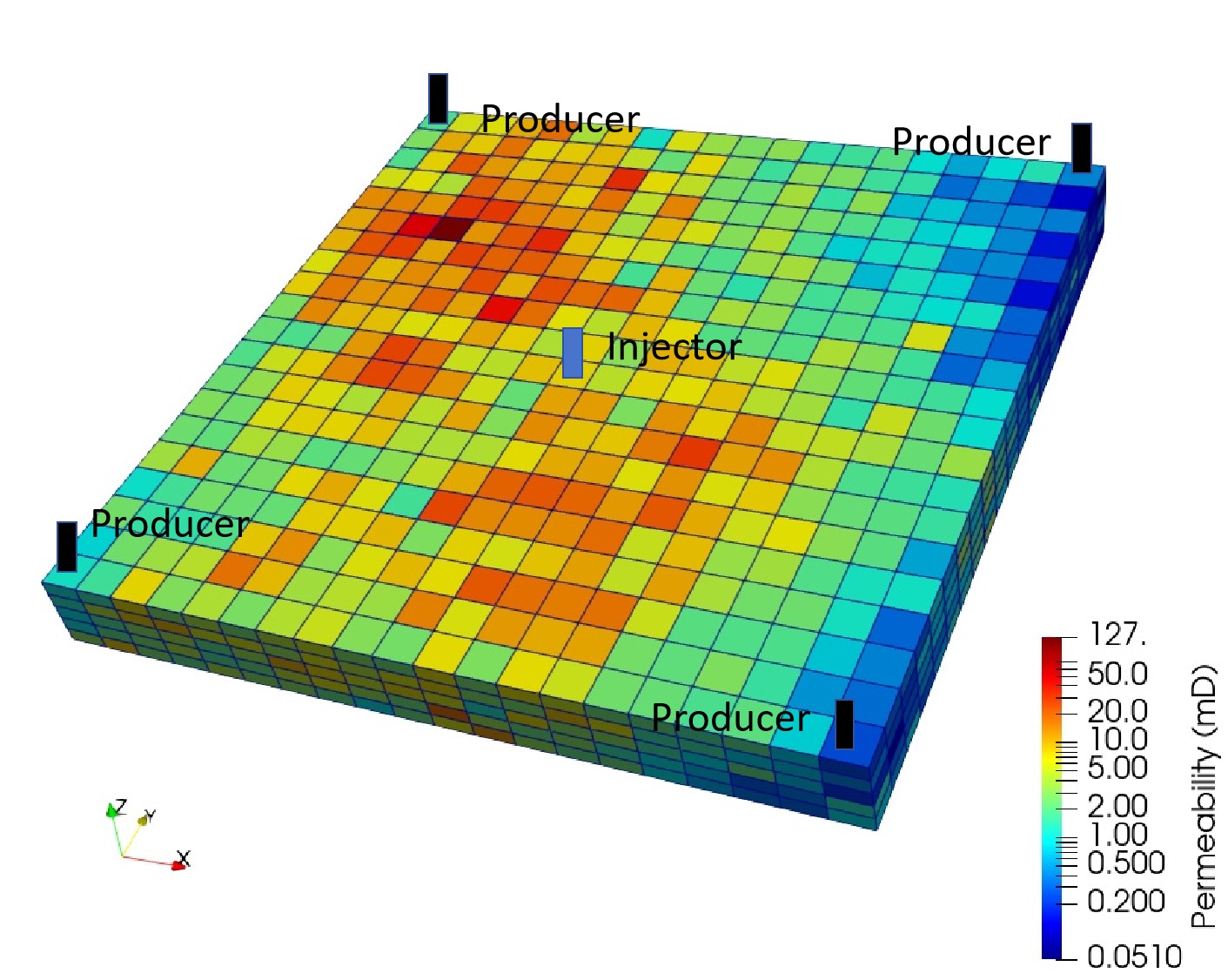}
	\caption{The ground truth for generating the observed flow rate data in test case 2.}
	\label{trueperm3d}
\end{figure}

\begin{figure}[!htp] 
	\centering
	\begin{subfigure}[t]{0.49\textwidth}
		\centering
		\includegraphics[width=1\textwidth]{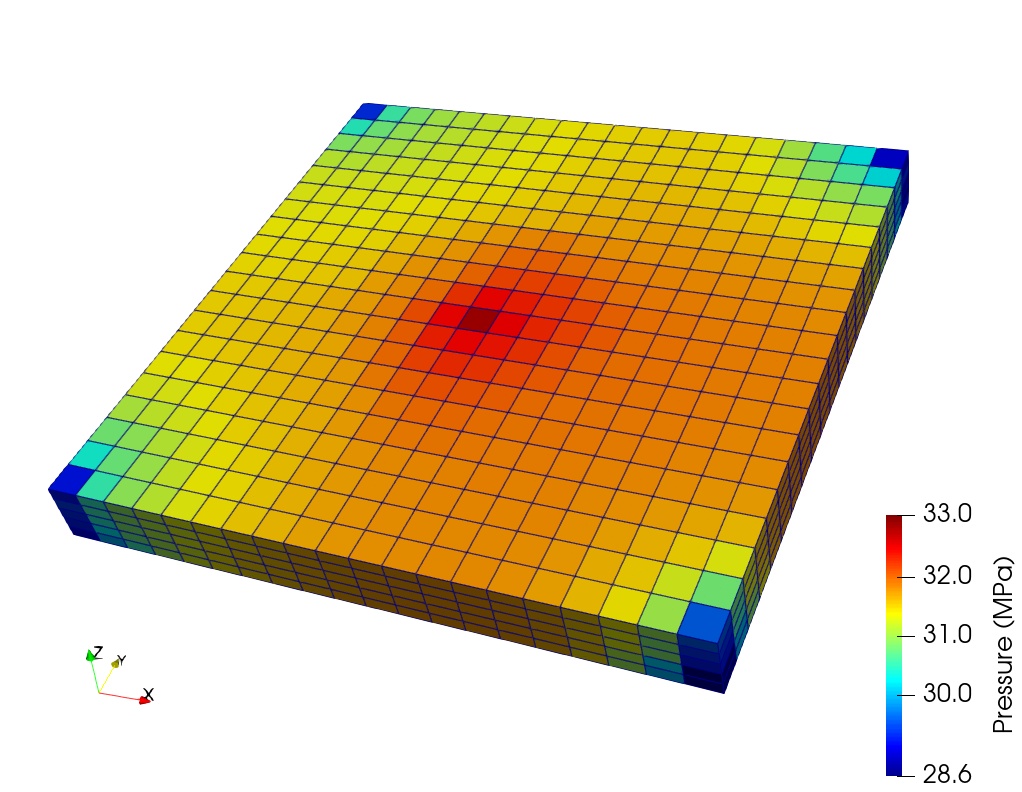}
		\caption{}
		\label{truepress3d} 
	\end{subfigure}
	\begin{subfigure}[t]{0.49\textwidth}
		\centering
		\includegraphics[width=1\textwidth]{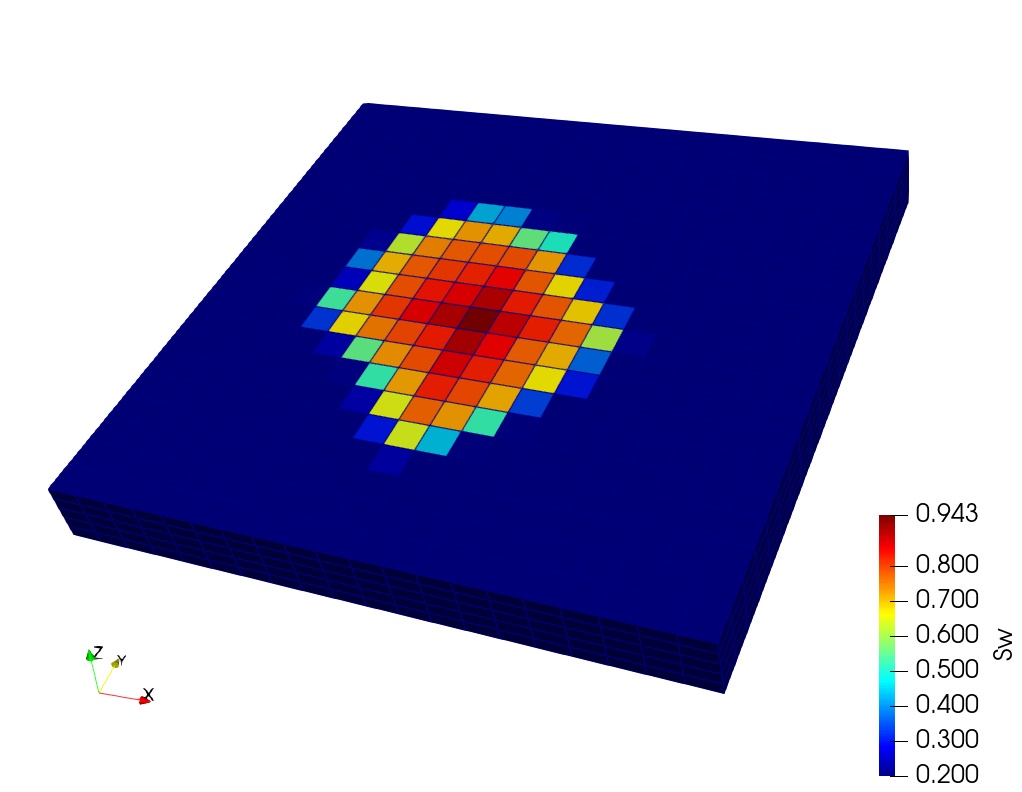}
		\caption{}
		\label{truesw3d}
	\end{subfigure}
	\caption{The pressure (b) and saturation (b) fields at end of simulation using the true permeability field for test case 2.}
\end{figure}

The injector in the centre of the reservoir is at constant flow rate and the four producers at corners are at constant bottom-hole pressure. The ground truth for the permeability field and the five wells are shown in Fig.~\ref{trueperm3d}. The pressure and saturation fields at end of simulation for the ground truth are shown in Figs.~\ref{truepress3d} and \ref{truesw3d}, respectively. For forward simulation, there are 2160 time steps with $dt=10k$ seconds. 
The permeability field is uncertain and all other constant parameters are shown in Table \ref{para2}. 

For NEI, we conduct forward simulations on the 100 prior realisations for once. Let each subset contain no more than four realisations, a total of $\sum_{i=1}^{4}C_{100}^i=4087975$ expectations for flow rates on subsets are computed. Next, the upper and lower expectations of the subset expectation are calculated to determine whether they fit the upper and lower bounds of the observed data. Based on this assessment, an appropriate threshold for loss is selected, resulting in 190 posterior subsets on which the expected flow rates cover the varying interval of observation. The posterior subsets contain 44 nonrepetitive realizations which are simulated further for prediction.

\begin{figure}[!htp] 
	\centering
	\begin{subfigure}[t]{0.48\textwidth}
		\centering
		\includegraphics[width=1\textwidth]{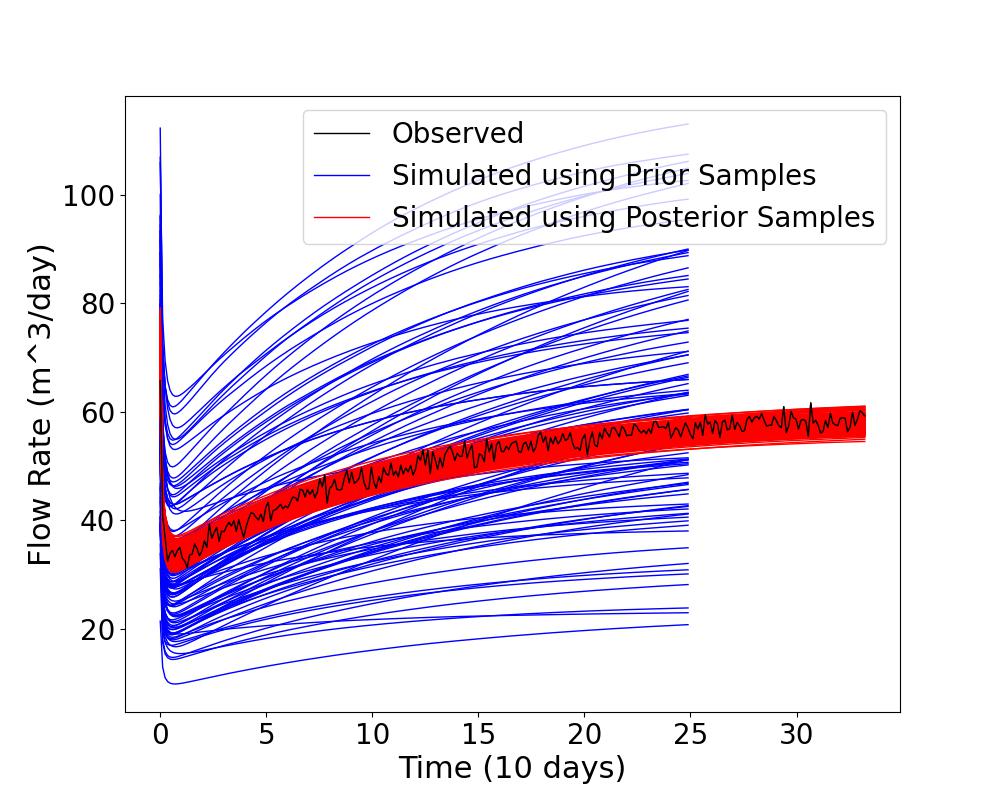}
		\caption{}
		\label{comparewell1} 
	\end{subfigure}
	\begin{subfigure}[t]{0.48\textwidth}
		\centering
		\includegraphics[width=1\textwidth]{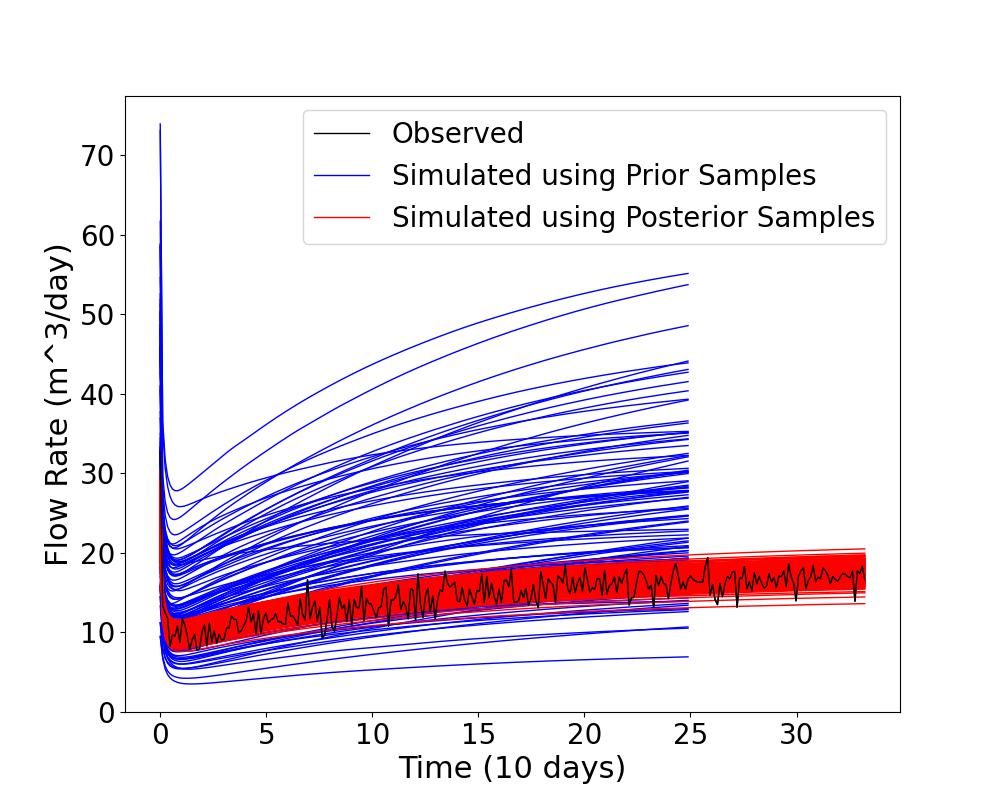}
		\caption{}
		\label{comparewell2}
	\end{subfigure}
	\begin{subfigure}[t]{0.48\textwidth}
		\centering
		\includegraphics[width=1\textwidth]{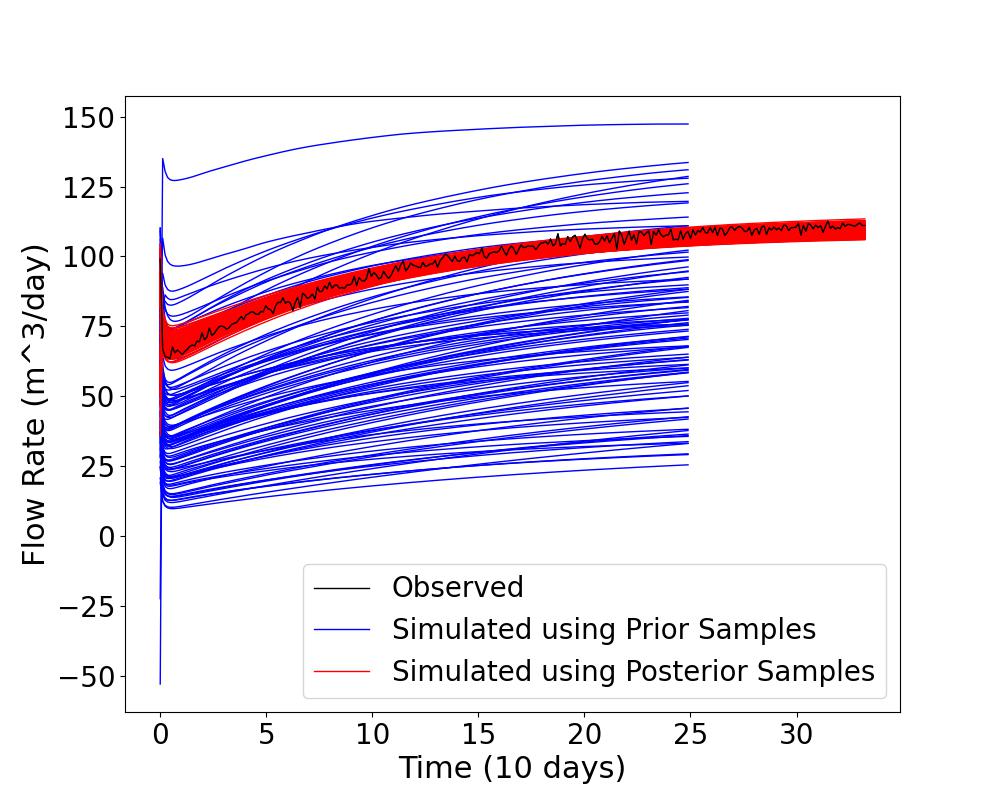}
		\caption{}
		\label{comparewell3}
	\end{subfigure}
	\begin{subfigure}[t]{0.48\textwidth}
		\centering
		\includegraphics[width=1\textwidth]{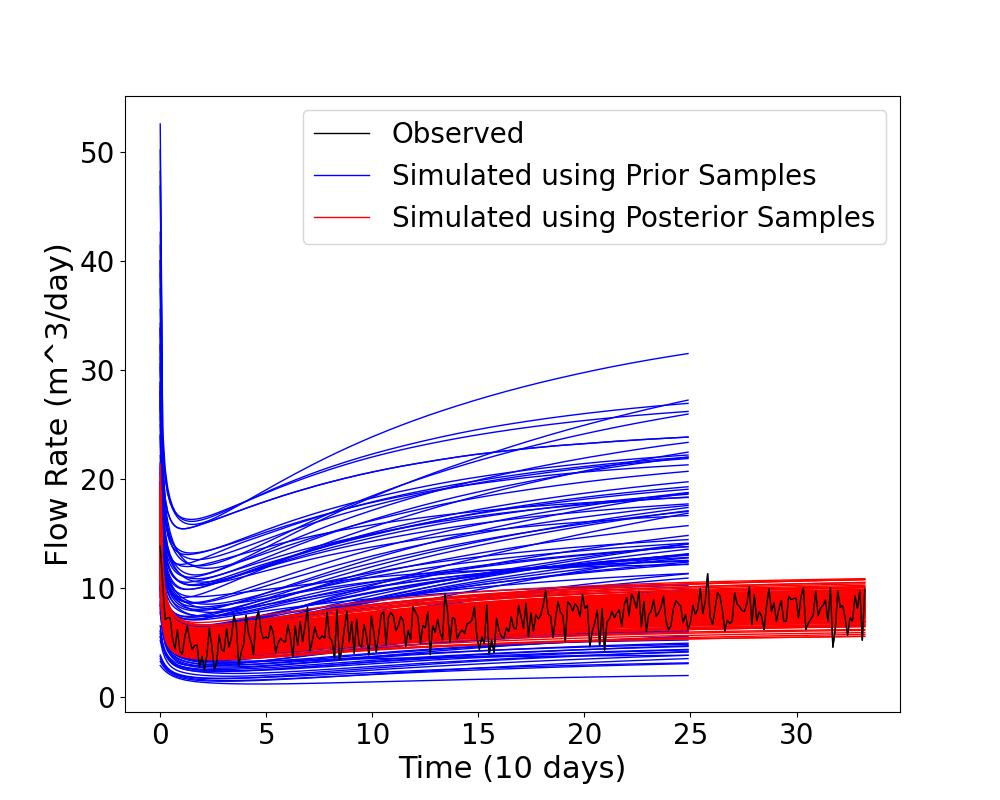}
		\caption{}
		\label{comparewell4}
	\end{subfigure}
	\caption{The simulated flow rates at producers using prior realisations and posterior subsets by NEI for test case 2.}
	\label{comparewells3D}
\end{figure}

\begin{figure}[!htp] 
	\centering
	\begin{subfigure}[t]{0.48\textwidth}
		\centering
		\includegraphics[width=1\textwidth]{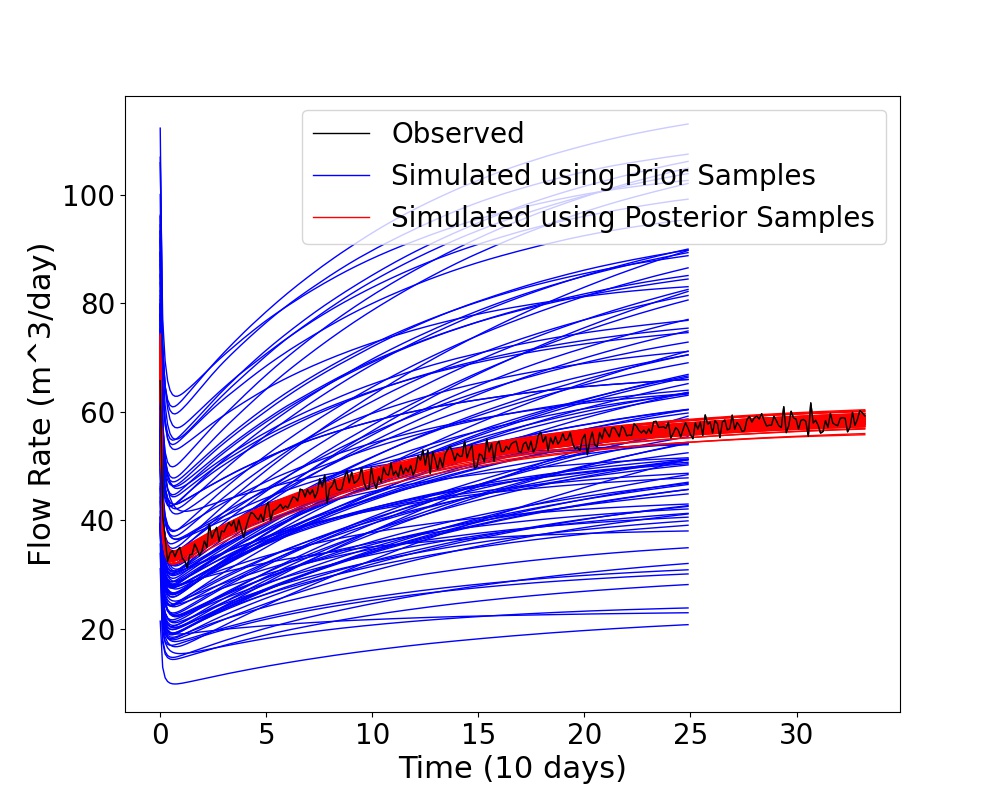}
		\caption{}
	\end{subfigure}
	\begin{subfigure}[t]{0.48\textwidth}
		\centering
		\includegraphics[width=1\textwidth]{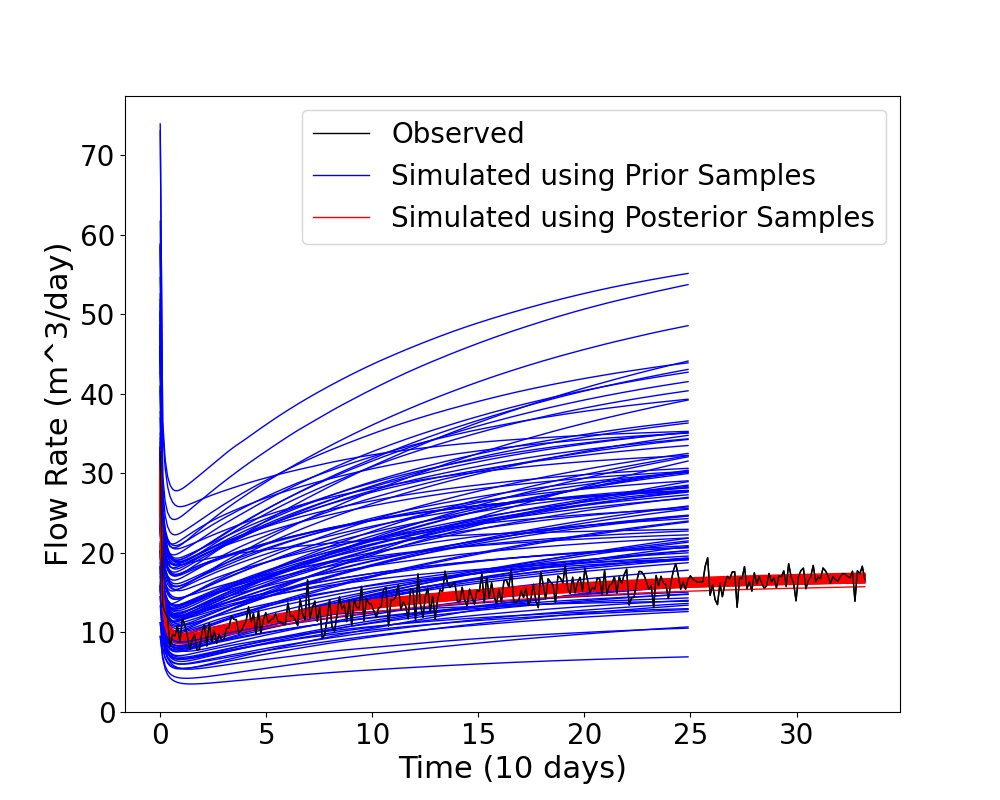}
		\caption{}
	\end{subfigure}
	\begin{subfigure}[t]{0.48\textwidth}
		\centering
		\includegraphics[width=1\textwidth]{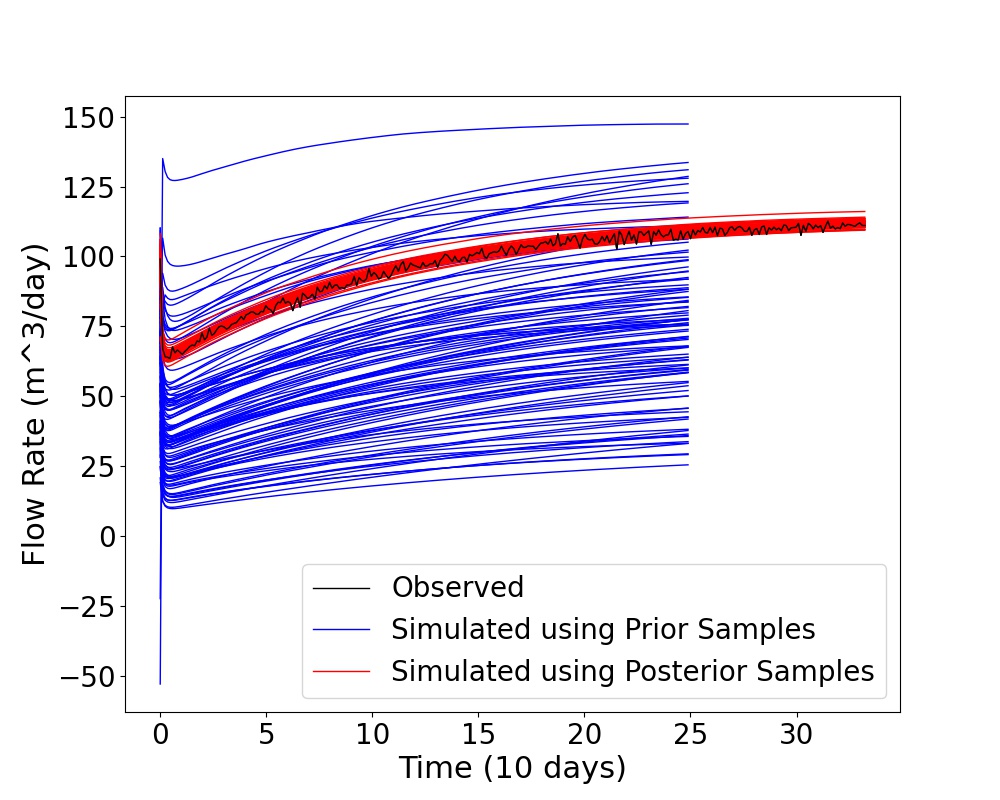}
		\caption{}
	\end{subfigure}
	\begin{subfigure}[t]{0.48\textwidth}
		\centering
		\includegraphics[width=1\textwidth]{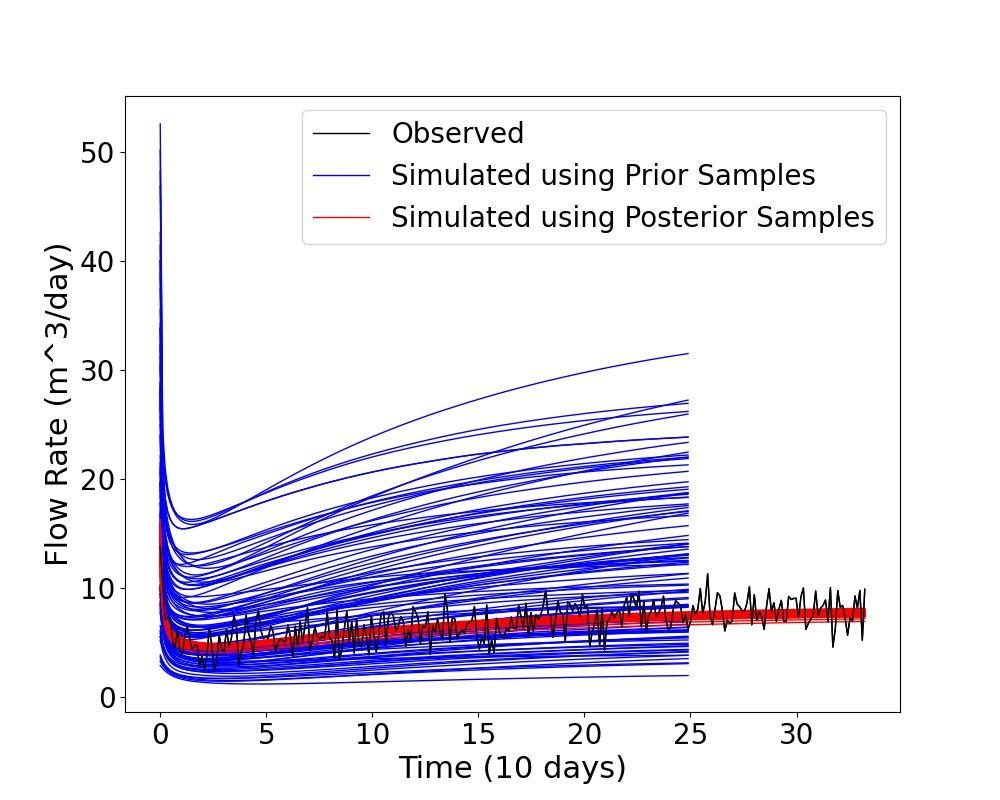}
		\caption{}
	\end{subfigure}
	\caption{The simulated flow rates at wells using prior and posterior realisations by ESMDA for test case 2.}
	\label{comparewells3Desmda}
\end{figure}

For ESMDA, the realizations are updated for two iterations, requiring an additional 200 forward model runs. The inference and predicted flow rates using NEI and ESMDA are shown in Figs.~\ref{comparewells3D} and \ref{comparewells3Desmda}, respectively.

The computational cost for NEI is about 6.45 hours, where 4.45 hours is for 100 forward simulation runs on prior realisations, and 2 hours is for evaluating the expectations on subsets (Table \ref{time}). The computational cost for ESMDA is about 13.5 hours for 300 forward simulation runs. In addition, the prediction efficiency of NEI is higher than ESMDA as forward simulation runs are conducted on less realisations. 

\begin{table}[!htp]
	\centering
	\begin{tabular}{|l|l|}
		\hline
		Porosity & 20\%\\
		\hline
		Water Compressibility & $4\times 10^{-6}$ Pa$^{-1}$\\
		\hline
		Oil Compressibility & $100\times 10^{-6}$ Pa$^{-1}$\\
		\hline
		Water Dynamic Viscosity & 0.001 Pa*s\\
		\hline
		Oil Dynamic Viscosity & 0.0018 Pa*s\\
		\hline
		Irreducible Water Saturation & 0.2\\
		\hline
		Production Index effective distance & 13.29 m\\
		\hline
		Bottom-hole Pressure & 28 MPa \\
		\hline
	\end{tabular}
	\caption{A summary of constant physical and well parameters for test case 2.}
	\label{para2}
\end{table}

\subsection{Test Case 3: Single-Phase Transient Flow in 3D Channelized Reservoir}
The third example is the inversion of the permeability field of EGG reservoir model \citep{jansen2014egg}. The EGG model is a synthetic channelized reservoir model consisting of 101 realisations considered to be geologically meaningful. The permeability field is non-Gaussian. Given the limited prior realisations, its probability distribution is effectively uncertain.

\begin{figure}[!htp] 
	\centering
	\begin{subfigure}[t]{0.48\textwidth}
		\centering
		\includegraphics[width=1\textwidth]{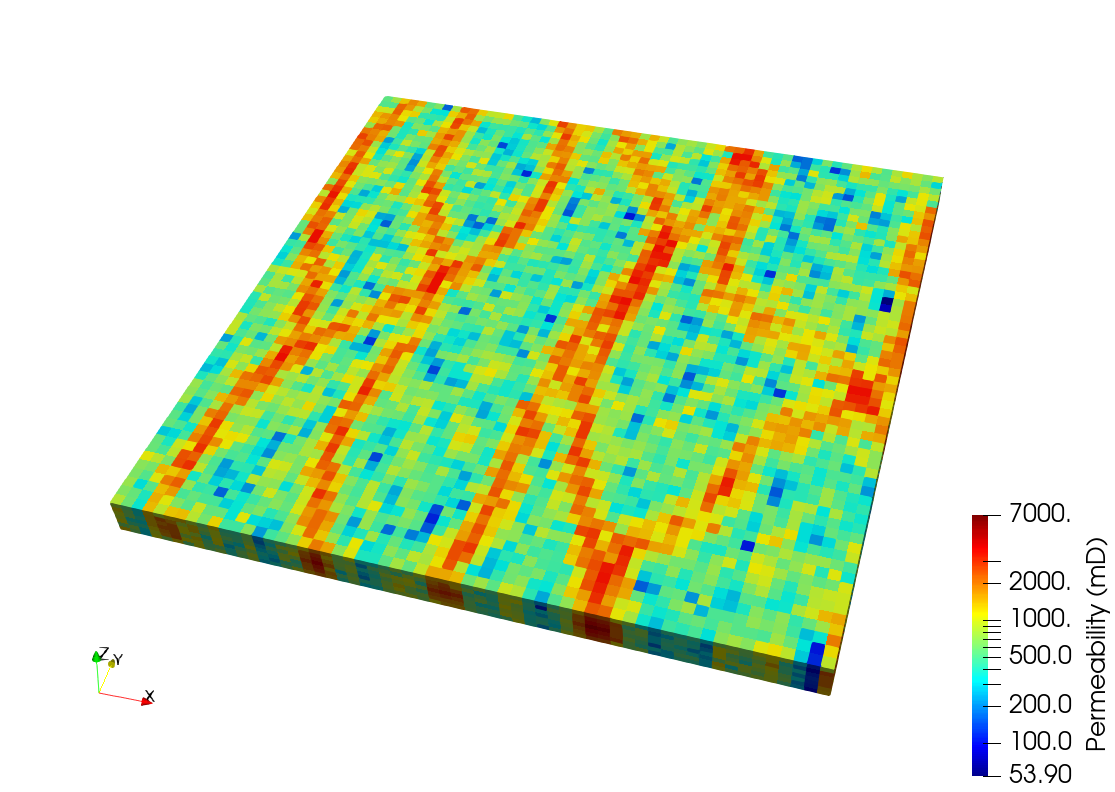}
		\caption{}
		\label{eggperm} 
	\end{subfigure}
	\begin{subfigure}[t]{0.48\textwidth}
		\centering
		\includegraphics[width=1\textwidth]{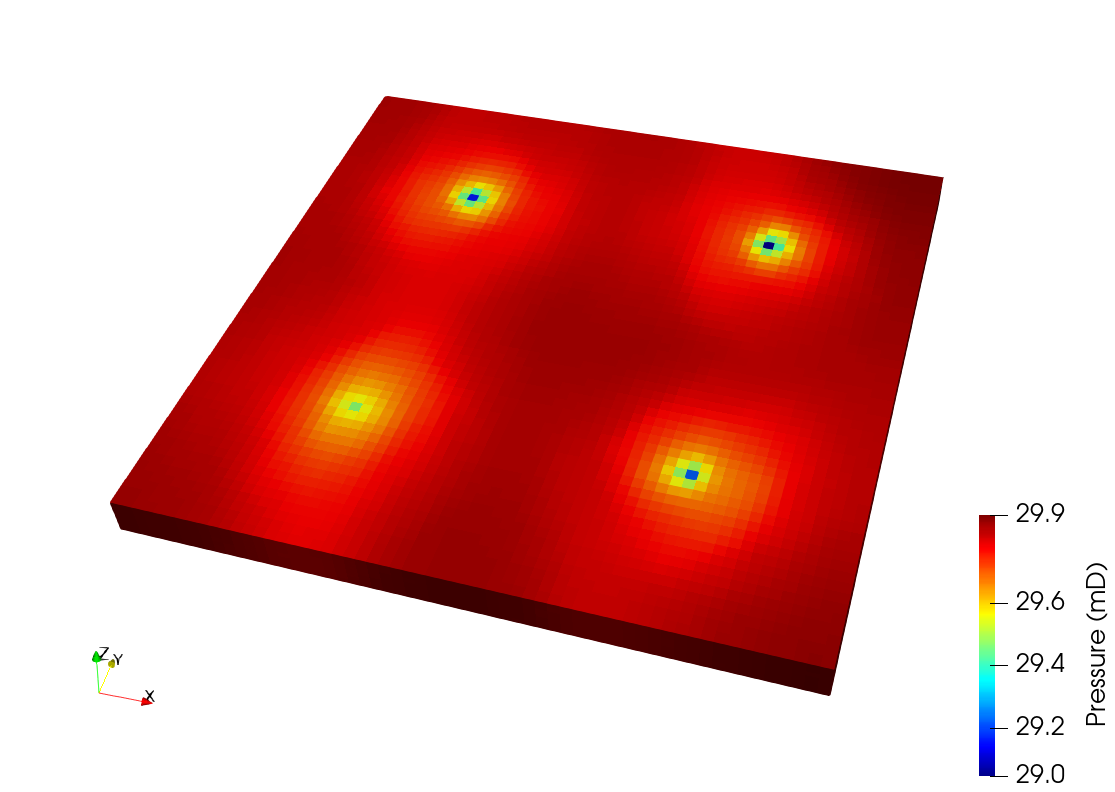}
		\caption{}
		\label{eggpress}
	\end{subfigure}
	\caption{The ground truth for the permeability field of the EGG reservoir model (a) and the pressure field at the end of simulation (b) in test case 3.}
	\label{eggtrue}
\end{figure}

\begin{figure}[!htp] 
	\centering
	\begin{subfigure}[t]{0.48\textwidth}
		\centering
		\includegraphics[width=1\textwidth]{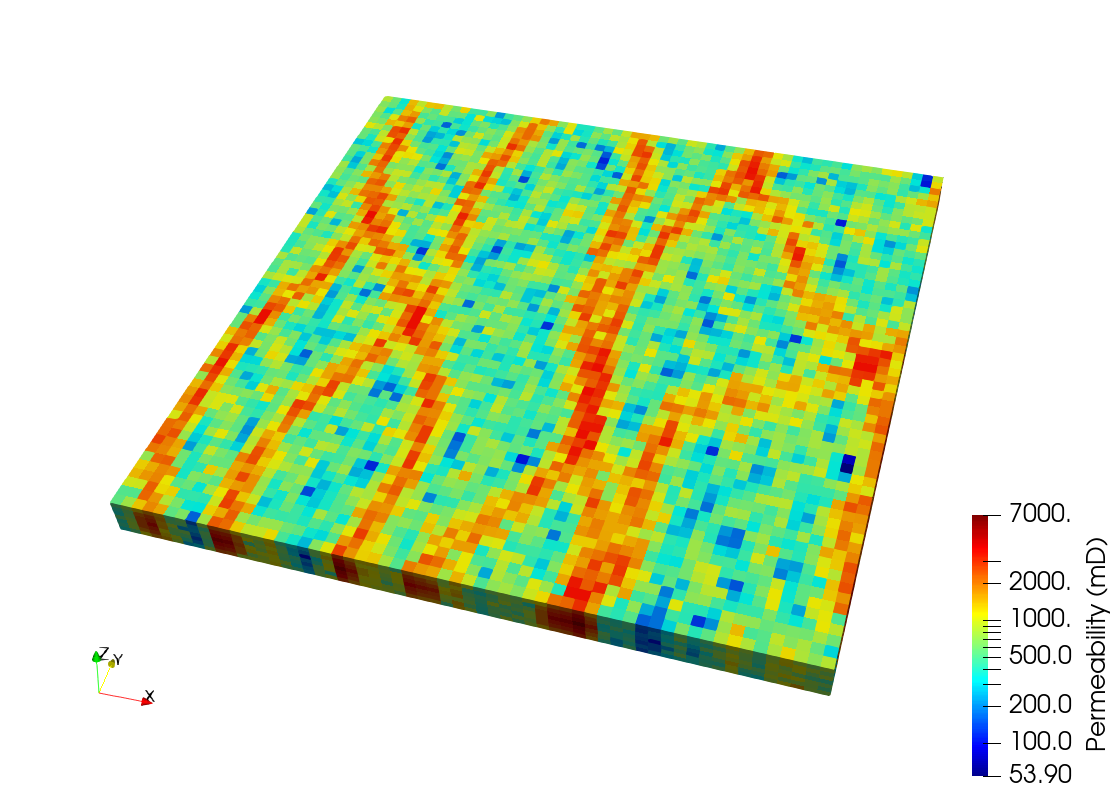}
		\caption{}
	\end{subfigure}
	\begin{subfigure}[t]{0.48\textwidth}
		\centering
		\includegraphics[width=1\textwidth]{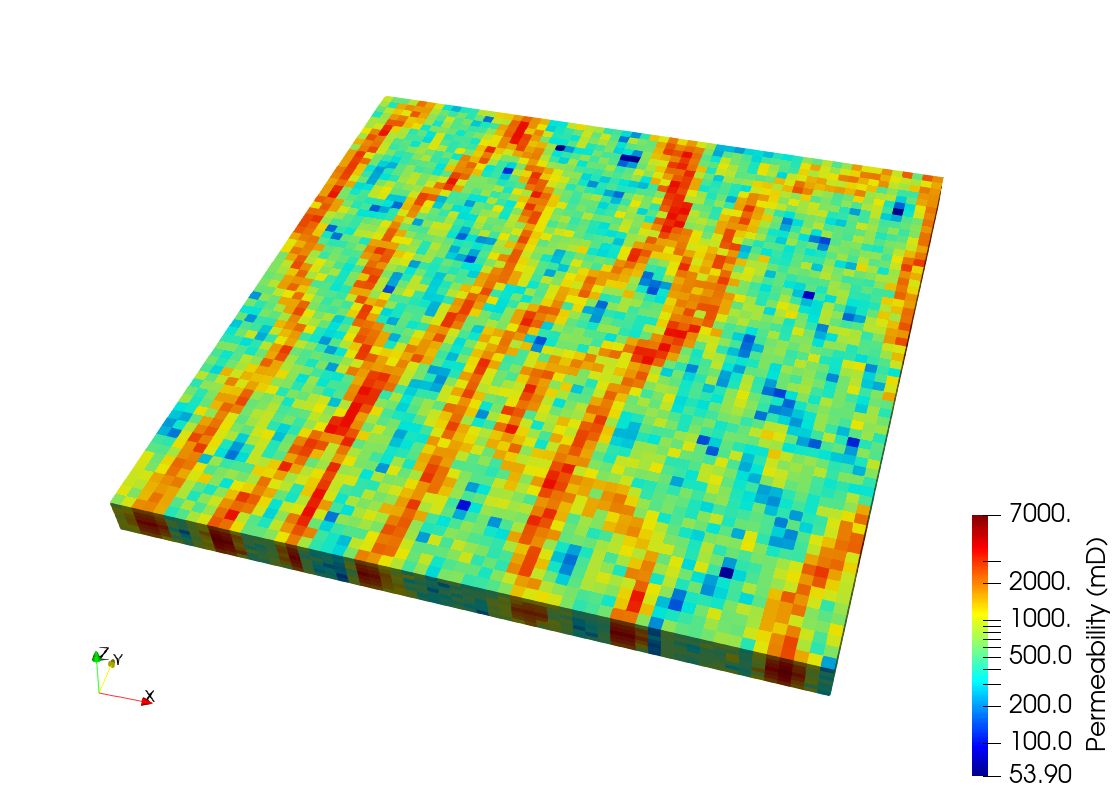}
		\caption{}
	\end{subfigure}
	\caption{Two of the prior realisations for the EGG reservoir model in test case 3.}
	\label{eggprior}
\end{figure}

One realisation is used as the ground truth, while all others serve as the prior realisations. The reservoir size is 480 m$\times$480 m$\times$28 m discretized by a 60$\times$60$\times$7 grid where all cells are active in the current study. The initial condition is that all cells are of pressure 30 MPa. There are four wells in the reservoir producing at the same bottom-hole pressure 26 MPa. The constant physical and well parameters are shown in Table \ref{para3}. The ground truth for the permeability field for generating the synthetic 'observed' flow rates with noise is shown in Fig.~\ref{eggperm}. The pressure field at the end of simulation is visualised in Fig.~\ref{eggpress}. For forward simulation, there are 72 time steps with $dt=5$k seconds. Two of the prior realisations are shown in Fig.~\ref{eggprior} where the main characteristics of the reservoir are kept but the channels are of different location and shapes from the ground truth. 

For NEI, we run forward simulations on the 100 prior models. Let each subset contain no more than three realizations, a total of $\sum_{i=1}^{3}C_{100}^i=166750$ subsets are built and evaluated to obtain the expectations of flow rates. Next, the upper and lower expectations are calculated to determine whether they fit the upper and lower bounds of the observed data. Based on this assessment, an appropriate threshold for loss is selected, resulting in 1385 posterior subsets with associated expectations covering the varying interval of observation. 

\begin{figure}[!htp] 
	\centering
	\begin{subfigure}[t]{0.48\textwidth}
		\centering
		\includegraphics[width=1\textwidth]{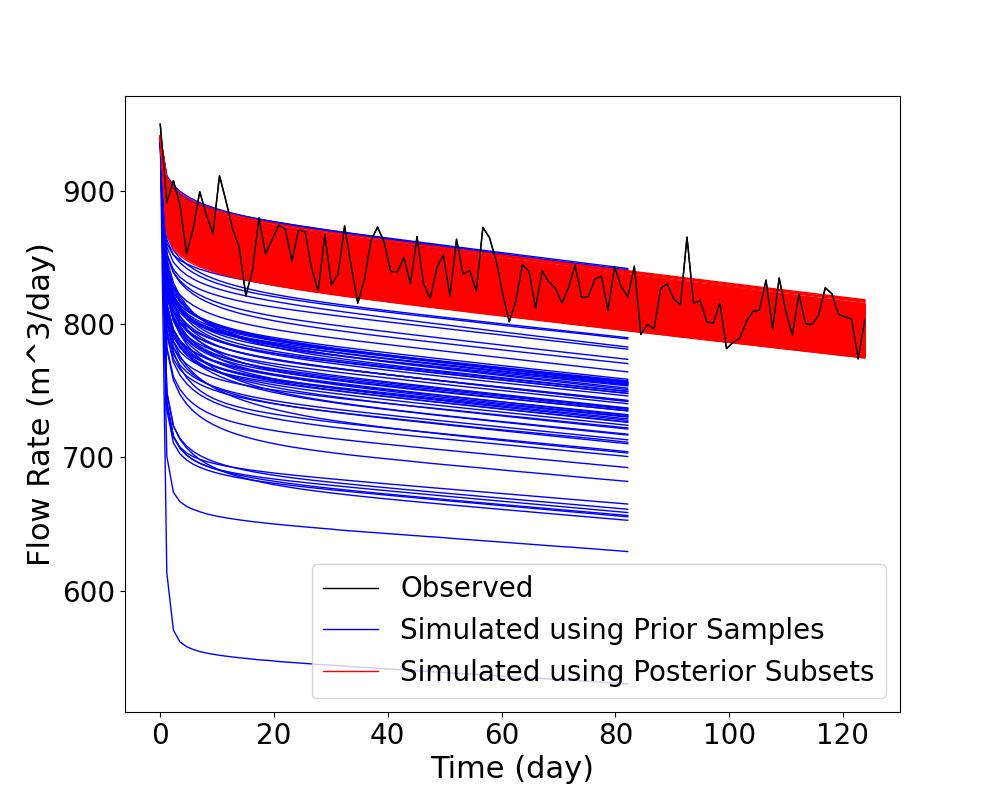}
		\caption{}
	\end{subfigure}
	\begin{subfigure}[t]{0.48\textwidth}
		\centering
		\includegraphics[width=1\textwidth]{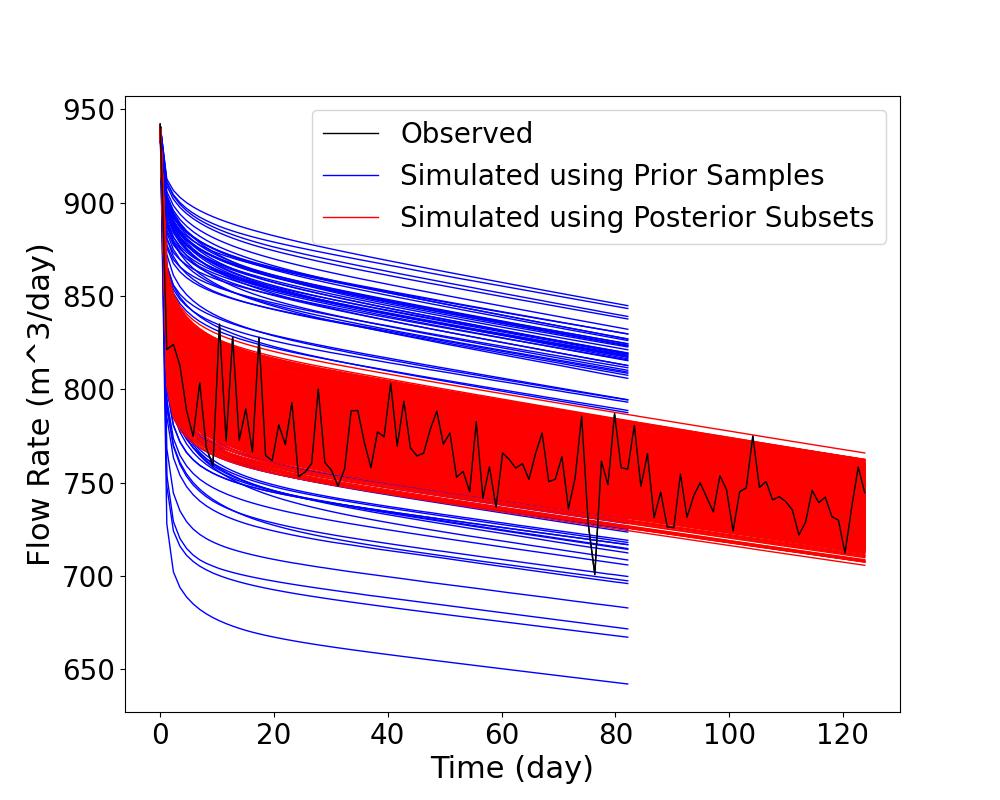}
		\caption{}
	\end{subfigure}
	\begin{subfigure}[t]{0.48\textwidth}
		\centering
		\includegraphics[width=1\textwidth]{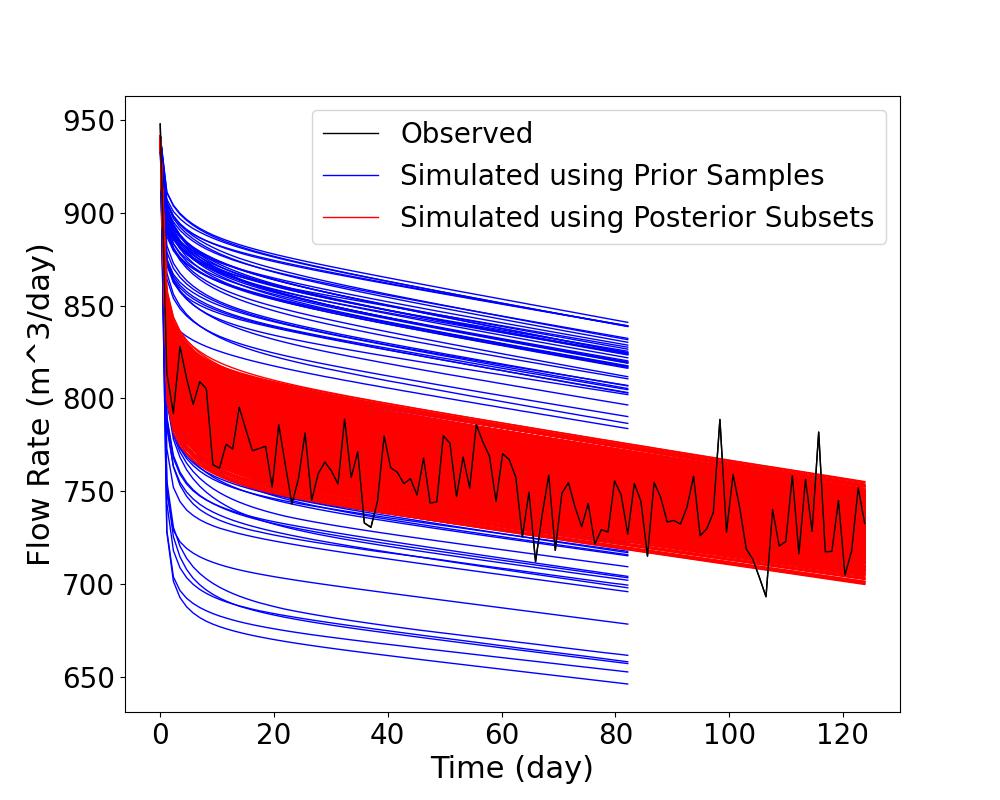}
		\caption{}
	\end{subfigure}
	\begin{subfigure}[t]{0.48\textwidth}
		\centering
		\includegraphics[width=1\textwidth]{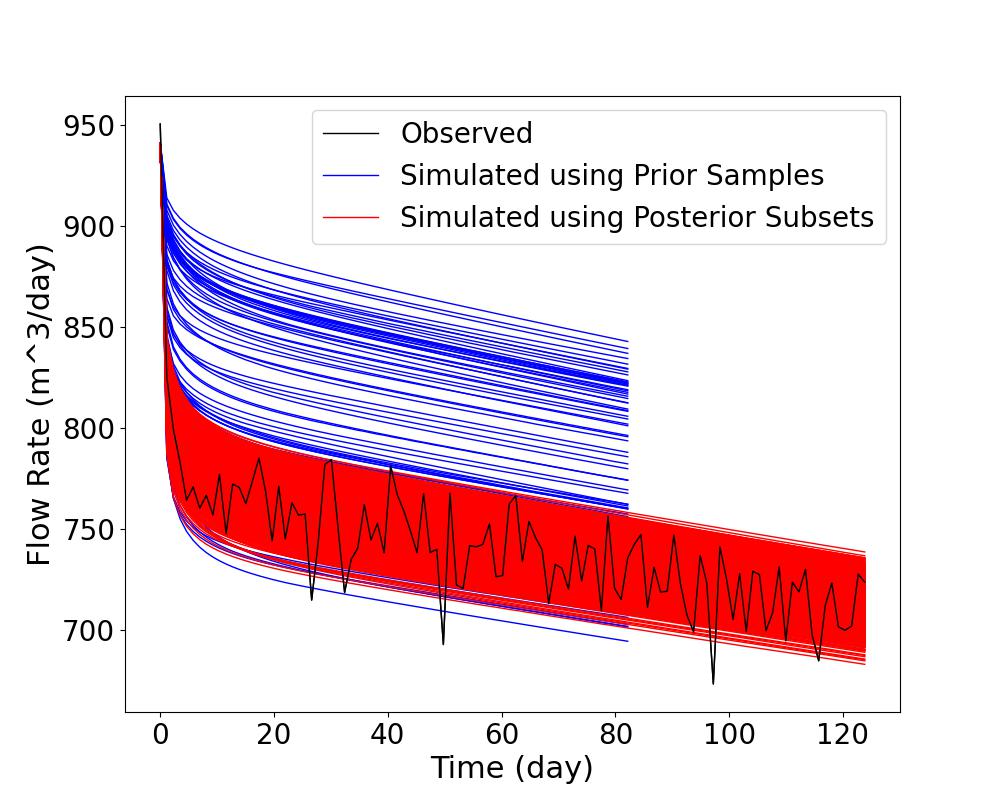}
		\caption{}
	\end{subfigure}
	\caption{The simulated flow rates at wells using prior realisations and posterior subsets by NEI for test case 3. There are up to three realisations in each subset.}
	\label{comparewellseggsize3}
\end{figure}

Flow rates simulated using prior realisations and posterior subsets are shown in Fig.~\ref{comparewellseggsize3}. Since there are 72 non-repetitive realisations in the union of the posterior subsets, we only need to simulate the 72 realisations to compute respective expectations for prediction. The computational cost of NEI is around 16 hours for simulating the 100 prior realisations plus 1 minute for evaluating expectations on the subsets (Table \ref{time}). The upper and lower expectations are the bounds of the predicted flow rates quantifying the uncertainty of prediction.

On the other hand, it is challenging for ESMDA to converge in the current study. 

\begin{table}[!htp]
	\centering
	\begin{tabular}{|l|l|}
		\hline
		Porosity & 20\%\\
		\hline
		Total Compressibility & $5\times 10^{-8}$ Pa$^{-1}$\\
		\hline
		Dynamic Viscosity & 0.002 Pa*s\\
		\hline
		Production Index & 3.891$\times$10$^{-4}$ m$^3$Sec$^{-1}$MPa$^{-1}$\\
		\hline
		Bottom-hole Pressure & 26 MPa \\
		\hline
	\end{tabular}
	\caption{A summary of constant physical and well parameters for test case 3.}
	\label{para3}
\end{table}

\newpage
\section{Conclusions}
In the current study, a new NEI method has been established for efficient UQ and history matching of subsurface transient Darcy flows accounting for distribution uncertainty for physical parameters. In NEI, no repetitive runs of the forward model is needed in the inference process. Forward runs are conducted on the prior realisations once followed by computation of expectations in the data space. Given observed data, the loss of each expectation is evaluated, then posterior expectations are obtained using a robust threshold to make sure that the uncertainty interval of observation is approximately covered by posterior expectations. The posterior uncertainty is therefore quantified by the nonlinear expectations defined as the limits of expectations. For prediction, a reduced number of non-repetitive realisations in the posterior subsets are simulated further and the uncertainty of prediction can be quantified using nonlinear expectations. The robustness and efficiency of NEI have been validated using test cases of of single- and two-phase transient Darcy flows in 2D and 3D heterogeneous reservoirs with Gaussian and non-Gaussian parameter random fields. Assuming that the observed data can be covered by the range of simulation results using prior realisations, NEI is not only highly efficient, but can converge on challenging history matching problems in channelized reservoirs with complex geological structures presenting distribution uncertainty.

\section*{Code Availability}
The NEI codes are available at "https://github.com/zhaozhang2022/nonlinearExp". The prior realisations can be generated using SGeMS or downloaded according to \citet{jansen2014egg}. The dynamic responses for prior realisations can be obtained using a numerical simulator for solving the governing PDE. 

\section*{Authorship Statement}
Zhao Zhang: Conceptualization, Funding acquisition, Formal analysis, Investigation, Methodology, Validation, Writing-original draft. Xinpeng Li: Formal analysis, Methodology, Validation, Writing–review \& editing. Menghan Li: Investigation, Methodology. Jiayu Zhai: Formal analysis. Piyang Liu: Investigation, Validation. Xia Yan: Validation. Kai Zhang: Supervision.

\section*{Acknowledgements}
The research is supported by the Natural Science Foundation of Shandong Province (No.ZR2024MA057), the Fundamental Research Funds for the Central Universities and the Future Plan for Young Scholars of Shandong University.
\newpage



  \bibliographystyle{elsarticle-harv} 
  \bibliography{mybib}



\newpage

\end{document}